\documentclass[11pt, a4paper]{article}

\usepackage{amsmath}\multlinegap=0pt
\usepackage{amssymb}
\usepackage{dsfont}
\usepackage[francais,english]{babel}

\providecommand\mathbb{\bf}

\newcommand\R{{\mathbb R}}

\newcommand\N{{\mathbb N}}

\renewcommand\S{{\mathbb S}}

\newcommand\pref[1]{(\ref{#1})}

\newcommand\indi{\mathds{1}}

\newcommand\qtext[1]{\quad\text{#1}\quad}

\let \longto\longrightarrow

\let \eps\varepsilon
\let \phi\varphi
\newtheorem{thm}{Theorem}

\newtheorem{prop}{Proposition}
\newtheorem{lemma}{Lemma}

\newcounter{Remark}
\newenvironment{Remark}
     {\par\medbreak \refstepcounter{Remark}
       \noindent\textbf{Remark~\arabic{Remark}.}}
      {\hbox{ }\medbreak}
\newenvironment{proof}[1][]
      {\par\medbreak{\noindent\bfseries Proof#1.\quad}}
      {\hbox{}\hfill\fbox{\ }\bigbreak}
\newcounter{steps}

\newcommand\MM{{\cal{M}}}

\newcommand\DD{({\cal{D}})}
\newcommand\DDr{({\cal{D}}_r)}
\newcommand\DDk{({\cal{D}}_k)}
\newcommand\DDrk{({\cal{D}}_{r,k})}
\newcommand\MKrk{({\cal{MK}}_{r,k})}

\newcommand\MK{({\cal{MK}})}
\newcommand\MKk{({\cal{MK}}_k)}

\newcommand\EE{{\cal E}}
\newcommand\LL{{\mathcal{L}}}

\newcommand\M{{\cal M}}
\newcommand\HH{{\cal H}}

\newcommand\idd{\int_{\R^{d\times d}}}

\newcommand\idpbd{\int_{\R^{d\times d}\setminus B_r^d}}

\newcommand\ovga{\overline{\gamma}}
\newcommand\ovH{\overline{H}}
\newcommand\undH{\underline{H}}

\newcommand\ogm{\overline{\gamma}}

\newcommand\leb[2]{\LL^#1_\text{#2}}
\def\<#1,#2>{\left<#1,#2\right>}
\let\bar\overline

\title {Optimal transportation for the determinant}
\author {G.~Carlier, B. Nazaret
\thanks{\scriptsize\ CEREMADE, UMR CNRS 7534, Universit\'e Paris Dauphine, Pl. de Lattre de Tassigny, 75775 Paris Cedex 16, FRANCE \texttt{carlier@ceremade.dauphine.fr}, \texttt{nazaret@ceremade.dauphine.fr}. }}
\date{4th december 2006}

\begin{document}

\maketitle

\begin{abstract}
Among $\R^3$-valued triples of random vectors $(X,Y,Z)$ having fixed marginal probability laws, what is the best way to jointly draw $(X,Y,Z)$ in such a way that the simplex generated by $(X,Y,Z)$ has maximal average volume?  Motivated by this simple question, we study optimal transportation problems with several marginals when the objective function is the determinant or its absolute value.
\end{abstract}



\section{Introduction}

Given two probability measures $\mu_1$ and $\mu_2$ on $\R^d$ and some objective function $H$ : $\R^d\times \R^d\rightarrow \R$, the classical Monge-Kantorovich optimal transportation problem consists in finding a probability measure $\gamma$ on $\R^d\times \R^d$ having $\mu_1$ and $\mu_2$ as marginals (i.e. a \emph{transportation plan} between $\mu_1$ and $\mu_2$) maximizing the total objective $\int_{\R^d\times \R^d} H(x,y)d\gamma(x,y)$. In his famous article \cite{bre}, Brenier solved the case $H(x,y)=\<x,y>$ and proved (under mild regularity assumptions) that there is a unique optimal transportation plan which is further characterized by the property of being supported by the graph of the gradient of some convex function. Brenier's seminal results have been extended to the case of more general costs (see McCann and Gangbo \cite{mcgan}) and the subject has received a lot of attention in the last 15 years because of its numerous applications in fluid mechanics, probability and statistics, PDE's, shape optimization, mathematical economics.... The literature on this very active field of research is too vast to give an exhaustive bibliography here, we rather refer to the books of Villani \cite{villani} and Rachev and R\"{u}schendorf \cite{raru} and the references therein.

\smallskip

In the present article, we are interested in an optimal transportation problem with several marginals. Given  $d$ probability measures $\mu_1,....,\mu_d$ on $\R^d$ and an objective function $H$ : $\R^{d\times d}\rightarrow \R$, the problem is to find a  probability measure $\gamma$ on $\R^{d\times d}$ having $\mu_1,...,\mu_d$ as marginals  maximizing $\int_{\R^{d\times d} }H(x_1,...,x_d)d\gamma(x_1,...,x_d)$. In contrast with the case of two marginals, there are few results on the optimal transportation problem when more than two marginals are involved with the exception of the article of Gangbo and  \'Swi\c ech \cite{gansw} who fully solved the case $H(x_1,...,x_d):=\sum_{i,j} \<x_i,x_j>$. For this particular problem, Gangbo and  \'Swi\c ech proved existence and uniqueness of an optimal transportation plan which is supported by the graph of some transport map.  In the present paper, we will pay attention here to different objective functions, namely: $H(x)=\det(x)$ or $H(x)=\vert \det(x)\vert$.

\smallskip

The choice of such functions of the determinant is motivated by the following simple question:  among random $\R^3$-valued vectors $X$,$Y$, and $Z$, with fixed marginal probability laws, what is the best way to draw jointly $(X,Y,Z)$ so that the simplex with vertices $(0,X,Y,Z)$ has maximal average volume? Denoting by $\mu_1$, $\mu_2$ and $\mu_3$ the (fixed) probability laws of the random vectors $(X,Y,Z)$, the previous problem amounts to maximize 
\[\int_{\R^{3\times 3}} \vert \det(x,y,z)\vert d\gamma(x,y,z)\]
among joint probability laws $\gamma$ having  $\mu_1$, $\mu_2$ and $\mu_3$ as marginals. In some cases, we will see that solving the problem above amounts to solve the simpler problem where $\vert \det \vert$ is replaced by $\det$ and we will study this case in dimension $d\geq 2$. If $d=2$, we have $\det(x,y)=\<x,Ry>$ (with $R$ the rotation of angle $-\pi/2$), so that, up to the change of variables $y'=Ry$, the optimal transportation problem with the determinant is a special case of the problem solved by Brenier \cite{bre}.

\smallskip

Let us define some notations.  In the sequel, given $X$ a locally compact separable metric  space, we denote by $\MM(X)$ (respectively $\MM^1_+(X)$) the set of Radon measures (respectively of Radon probability measures) on $X$. If $X$ and $Y$ are locally compact separable metric spaces, $\mu\in \M_+^1(X)$,  and $f$ : $X\rightarrow Y$ is a Borel map we shall denote by $f\sharp \mu$ the push forward of $\mu$ through $f$ i.e. the element of $\M_+^1(Y)$ defined by $f\sharp \mu (B)=\mu (f^{-1}(B))$ for every Borel subset $B$ of $Y$.  If $\gamma\in \MM(\R^{d\times d})$ and $i\in\{1,...,d\}$,  $\pi_i \sharp \gamma\in \MM(\R^d)$ is called the $i$-th \emph{marginal} of $\gamma$ (where $\pi_i$ is the $i$-th canonical projection). One can also define $\pi_i \sharp\gamma$ by:
\[\int_{\R^d} f(x_i) d (\pi_i \sharp \gamma)(x_i)=\int_{\R^{d\times d}} f(x_i)d\gamma(x_1,...,x_d),\]
for every bounded and continuous  function $f$ on $\R^d$. Given $d$ probability measures on $\R^d$, $\mu_1,...,\mu_d$ , we denote by $\Pi(\mu_1,...,\mu_d)$ the set of probability measures on $\R^{d\times d}$ having $\mu_1,...,\mu_d$ as marginals. In other words, $\gamma\in \MM^1_+(\R^{d\times d})$ belongs to $\Pi(\mu_1,...,\mu_d)$ if and only if
\[\int_{\R^{d\times d}} f(x_i)d\gamma(x_1,...,x_d)=\int_{\R^{d}} f(x_i)d\mu_i(x_i),\; \forall i=1,...,d,\]
for every bounded and continuous  function $f$ on $\R^d$. Given $H\in C^0(\R^{d\times d},\R)$, the Monge-Kantorovich optimal transportation problem with marginals $\mu_1,...,\mu_d$ and objective function $H$ then reads as:
\[ \sup_{\gamma \in \Pi(\mu_1,...,\mu_d)} \int_{\R^{d\times d}} H(x_1,...,x_d) d\gamma(x_1,...,x_d).\]
We shall focus here on the two special cases:
\[\MK  \sup_{\gamma \in \Pi(\mu_1,...,\mu_d)} \int_{\R^{d\times d}} \det(x_1,...,x_d) d\gamma(x_1,...,x_d)\]
and
\[ \MK_a \sup_{\gamma \in \Pi(\mu_1,...,\mu_d)} \int_{\R^{d\times d}} \vert \det (x_1,...,x_d) \vert d\gamma(x_1,...,x_d).\]

The paper is organized as follows. In section 2, we solve a particular example which actually gives insight on the general case. Section 3 is devoted to existence, duality and characterization of minimizers for $\MK$. In section 4, we construct minimizers in the case of radially symmetric marginals. In section 5, we give conditions ensuring that $\MK$ and $\MK_a$ are in fact equivalent. Section 6 contains various remarks regarding uniqueness issues.

\section{An elementary example}\label{uniform}

In this section, we study a simple but illustrative example, which actually contains most of the ideas necessary for the understanding of the general case. Let us consider the problem $\MK$ in dimension $3$, with $\mu_1=\mu_2=\mu_3=\leb{3}{B}$, where $\leb{3}{B}$ stands for the uniform probability measure on the unit ball $B$ of $\R^3$,

\begin{equation}
\label{3d-primal}
\sup_{\gamma\in\Pi(\leb{3}{B},\leb{3}{B},\leb{3}{B})}\left(\int_{B^3}\det(x,y,z)d\gamma(x,y,z)\right).
\end{equation}

Then, the following result holds.

\begin{thm}\label{thm-3d}
The problem \eqref{3d-primal} admits a solution $\bar\gamma$ given by, for every continuous $f:B^3\longto\R$,
\begin{equation}\label{3d-solution}
\int_{B^3}fd\bar\gamma=\frac{1}{|B|}\int_B\left(\int_{\S(x)}f\left(x,|x|y,x\wedge y\right)\frac{d\HH^1(y)}{2\pi}\right)dx,
\end{equation}
where
$$ \S(x)=\left\{ y\in \S^2 \qtext{s.t.} \<x,y>=0 \right\}.$$
\end{thm}

\begin{Remark}\label{rkor}
Let us remark that, the support of the optimal measure $\bar\gamma$ is the set of all the triples $(x,y,z)\in B^3$, such that $|x|=|y|=|z|$ and $(x,y,z)$ is a direct orthogonal basis of $\R^3$. As we shall see later on, this last property comes from the fact that the measures are radially symmetric. In addition, the definition of $\bar\gamma$  through its successive disintegrations  has the following probabilistic interpretation in terms of conditional laws. Consider $\bar\gamma$ as the law of a vector $(X,Y,Z)$ of random vectors in $\R^3$, each of them following an uniform law on the unit ball. Then, \eqref{3d-solution} means that the conditional probability of $Y$ given $X$ is uniform on ${S(X)}$ and that the  conditional probability of $Z$ given $(X,Y)$ is the Dirac mass at $\frac{X\wedge Y}{|X|}$.
\end{Remark}

\begin{Remark}\label{valeurabs}
Let us say a word about the case where the objective function is  $H(x,y,z)=|\det(x,y,z)|$ (volume maximization). For this $H$, $\bar\gamma$ is still a maximizer, but we could have chosen as well
$$
Z=-\frac{X\wedge Y}{|X|},
$$
for the third variable, and the solution is not unique. Actually, it is not the only source of non uniqueness, and we shall discuss this point in the last section.
\end{Remark}

\begin{proof}[ of theorem \ref{thm-3d}]

Assuming that  $\bar\gamma$ is admissible in \eqref{3d-primal}, we only need to prove that it is optimal. To do so, consider the variational problem
\begin{equation}\label{3d-dual}
\inf_{(\phi,\psi,\chi)\in\EE}\frac{1}{|B|}\left(\int_B\phi(x)dx+\int_B\psi(y)dy+\int_B\chi(z)dz\right),
\end{equation}
with $\EE$ the set of triples $(\phi,\psi,\chi)\in C^0(B,\R)^3$ such that
$$
 \phi(x)+\psi(y)+\chi(z)\geq\det(x,y,z),\; \forall (x,y,z)\in B^3.
$$
As we shall see in next section, \eqref{3d-dual} is dual to \eqref{3d-primal} in some sense. For all $(\phi,\psi,\chi)\in\EE$, and $\gamma\in\Pi(\leb{3}{B},\leb{3}{B},\leb{3}{B})$,
\begin{eqnarray*}
\int_{B^3}\det(x,y,z)d\gamma(x,y,z) & \leq & \int_{B^3}\left(\phi(x)+\psi(y)+\chi(z)\right)d\gamma(x,y,z),\\
& \leq & \frac{1}{|B|}\left(\int_B\phi(x)dx+\int_B\psi(y)dy+\int_B\chi(z)dz\right),
\end{eqnarray*}
and it immediately follows that
\begin{eqnarray}
\lefteqn{\sup_{\Pi(\leb{3}{B},\leb{3}{B},\leb{3}{B})}\int_{B^3}\det(x,y,z)d\gamma(x,y,z)} & & \nonumber\\
\label{dp-ineq} & & \quad\leq\inf_{\EE}\frac{1}{|B|}\left(\int_B\phi(x)dx+\int_B\psi(y)dy+\int_B\chi(z)dz\right).
\end{eqnarray}
Now, consider $(\phi_0,\psi_0,\chi_0)$ given by
$$
\forall x\in B, \ \phi_0(x)=\psi_0(x)=\chi_0(z)=\frac{|x|^3}{3}.
$$
Thanks to the Young inequality, we have
\begin{equation}\label{3d-young}
\forall (x,y,z)\in B^3, \ \phi_0(x)+\psi_0(y)+\chi_0(z) \geq |x||y||z| \geq \det(x,y,z),
\end{equation}
so that $(\phi_0,\psi_0,\chi_0)\in\EE$. In addition, according to remark \ref{rkor}, if $(x,y,z)\in\text{supp}(\bar\gamma)$, then all the inequalities in \eqref{3d-young} become equalities and then,
\begin{equation}\label{3d-condopt}
\forall (x,y,z)\in\text{supp}(\bar\gamma), \ \phi_0(x)+\psi_0(y)+\chi_0(z)=\det(x,y,z).
\end{equation}
It implies that \eqref{dp-ineq} is actually an equality and $\bar\gamma$ is optimal in \eqref{3d-primal}. To complete the proof, it just remains to show that $\bar\gamma$ is admissible, which is done in proposition \ref{3d-adm}.
\end{proof}

\begin{prop}\label{3d-adm}
Let $\bar\gamma$ defined as in theorem \ref{thm-3d}. Then, $\bar\gamma\in\Pi(\leb{3}{B},\leb{3}{B},\leb{3}{B})$.
\end{prop}
\begin{proof}
The fact that $\pi_1\sharp\bar\gamma=\leb{3}{B}$ is obvious by \eqref{3d-solution}. Now, let us suppose that we proved $\pi_2\sharp\bar\gamma=\leb{3}{B}$. Then, we can easily deduce the result for the third marginal. Indeed, since for every fixed $x\in B$, the map $y\mapsto\frac{x\wedge y}{|x|}$ is one-to-one from $\S(x)$ to itself and is simply a rotation with angle $\pi/2$, we have that for all continuous $f:B\to\R$,
\begin{eqnarray*}
\int_Bfd(\pi_3\sharp\bar\gamma) & = & \frac{1}{|B|}\int_B\left(\int_{\S(x)}f\left(x\wedge y\right)\frac{d\HH^1(y)}{2\pi}\right)dx,\\
& = & \frac{1}{|B|}\int_B\left(\int_{\S(x)}f(|x|y)\frac{d\HH^1(y)}{2\pi}\right)dx,\\
& = & \int_Bfd(\pi_2\sharp\bar\gamma).
\end{eqnarray*}
We then prove that $\pi_2\sharp\bar\gamma=\leb{3}{B}$. Let $f:B\to\R$ be a continuous function. Then
\begin{eqnarray*}
\int_Bfd(\pi_2\sharp\bar\gamma) & = & \frac{1}{|B|}\int_B\left(\int_{\S(x)}f(|x|y)\frac{d\HH^1(y)}{2\pi}\right)dx\\
& = & \frac{1}{|B|}\int_0^1r^2\left(\int_{\S^2}\left(\int_{\S(x)}f(ry)\frac{d\HH^1(y)}{2\pi}\right)d\HH^2(\sigma)\right)dr.\\
\end{eqnarray*}
Applying lemma \ref{fubsphere} proved in Appendix A, we get
\begin{eqnarray*}
\int_Bfd(\pi_2\sharp\bar\gamma) & = & \frac{1}{|B|}\int_0^1\left(\int_{\S^2}r^2f(ry)\left(\int_{S(y)}\frac{d\HH^1(\sigma)}{2\pi}\right)d\HH^2(y)\right)dr \\
& = & \int_B fd\leb{3}{B},
\end{eqnarray*}
which ends the proof.
\end{proof}

Let us notice that the condition \eqref{3d-condopt} completely characterizes the solutions of the dual problem \eqref{3d-dual}. In fact (see subsection \ref{cara}), up to the addition of constants that sum to $0$, $(\phi_0,\psi_0,\chi_0)$ is the unique solution of \eqref{3d-dual}. Yet, there are infinitely many solutions to \eqref{3d-primal}, indeed any $\gamma\in \Pi(\leb{3}{B},\leb{3}{B},\leb{3}{B})$ having its support in  all the triples $(x,y,z)\in B^3$, such that $|x|=|y|=|z|$ and $(x,y,z)$ is a direct orthogonal basis of $\R^3$ is optimal for \eqref{3d-dual} and we claim that there are infinitely many such probability measures (see section \ref{comments}).

\section{Duality, existence and characterization}

Given $d$ probability measures on $\R^d$, $\mu_1,...,\mu_d$ , we consider the problem 
\[\MK \;  \sup_{\gamma \in \Pi(\mu_1,...,\mu_d)} \int_{\R^{d\times d}} \det (x_1,...,x_d) d\gamma(x_1,...,x_d). \]
In this case, a natural assumption on the marginals  $\mu_1,...,\mu_d$, is the existence of $(p_1,...,p_d)\in (1,+\infty)^d$ such that:
\begin{equation}\label{hypintegr}
\sum_{i=1}^{d} \frac{1}{p_i}=1, \mbox{ and } \; a:=\sum_{i=1}^{d} \int_{\R^d} \frac{\vert x \vert^{p_i}}{p_i}d\mu_i(x)<+\infty.
\end{equation}
Defining for all $x=(x_1,....,x_d)\in \R^{d\times d}$:
\[H_0(x):=\sum_{i=1}^{d}  \frac{\vert x_i \vert^{p_i}}{p_i}\]
and using the fact  that $\vert \det(x)\vert \leq H_0(x)$, we immediately get 
\[-a\leq \int_{\R^{d\times d}} \det (x_1,...,x_d) d\gamma(x_1,...,x_d)\leq a,\; \forall \gamma\in \Pi(\mu_1,...,\mu_d).\]
Similarly, defining for all $x=(x_1,....,x_d)\in \R^{d\times d}$:
\begin{equation}\label{defhlower}
\ovH(x):=\det(x)+H_0(x),\; \undH(x):=\det(x)-H_0(x),
\end{equation}
we remark that for all $\gamma\in \Pi(\mu_1,...,\mu_d)$
\[\int_{\R^{d\times d}} \det (x) d\gamma(x)=\int_{\R^{d\times d}} \ovH d\gamma-a=\int_{\R^{d\times d}} \undH d\gamma+a.\]
The previous remark implies that when $H=\det$, one may without loss of generality replace $H=\det$ with $\ovH\geq 0$ or $\undH\leq 0$ in $\MK$.

\smallskip

A key point in Monge-Kantorovich theory, is to remark that (in a sense that will be made precise later), $\MK$ is dual to :
\[\DD \; \inf_{(\phi_1,...,\phi_d)\in \EE} \sum_{i=1}^d \int_{\R^d} \phi_i d\mu_i \] 
where $\EE$ is the set of $d$-uples of lower-semi continuous functions $(\phi_1,...,\phi_d)$ from $\R^d$ to $\R\cup \{+\infty\}$ such that
\begin{equation}\label{contraintedudual}
 \sum_{i=1}^d \phi_i(x_i)\geq \det(x_1,...,x_d),\; \forall \; (x_1,...,x_d)\in\R^{d\times d}. 
\end{equation}
Let us note that one obviously has $\inf \DD\geq \sup \MK$.

\subsection{The compact case}

In this paragraph, for further use, we consider the case of compactly supported marginals and of an arbitrary continuous objective function $H$. Denoting by $B_r$  the closed ball in $\R^d$, with center $0$ and radius $r$, we assume that $\mu_1,....,\mu_d$ are supported in  $B_r$ for some given $r>0$ and that $H$ is an arbitrary continuous function on $B_r^d$.  We then consider the optimal transportation problem:
\begin{equation}\label{mkh}
  \sup_{\gamma \in \Pi(\mu_1,...,\mu_d)} \int_{\R^{d\times d}} H(x_1,...,x_d) d\gamma(x_1,...,x_d). \end{equation}

We define its dual  by:
\begin{equation}\label{dualhr}
 \inf_{(\phi_1,...,\phi_d)\in \EE_r} \sum_{i=1}^d \int_{B_r} \phi_i d\mu_i \end{equation}
where $\EE_r$ is the set of $d$-uples of continuous functions $(\phi_1,...,\phi_d)$ from $B_r$ to $\R$ such that
\begin{equation}\label{contraintedudualr}
 \sum_{i=1}^d \phi_i(x_i)\geq H(x_1,...,x_d),\; \forall \; (x_1,...,x_d)\in B_r^d. 
\end{equation}

\begin{prop}\label{cascompact}
Assume that  $\mu_1,....,\mu_d$ are supported in the ball  $B_r$ and that $H$ is continuous on $B_r^d$  then both \pref{mkh} and \pref{dualhr} admit solutions and
\[\max \; \pref{mkh}= \min \; \pref{dualhr}.\]
Moreover, \pref{dualhr} admits a solution $(\phi_1,....,\phi_d)$ such that  for all $i=1,...,d$ and all $x_i\in B_r$, one has:
\begin{equation}\label{conj}
\phi_i(x_i)=\sup_{(x_j)_{j\neq i}\in B_r^{d-1}} \left\{H(x_1,...,x_d)-\sum_{j\neq i} \phi_j(x_j) \right\} 
\end{equation}
and, for all $x\in B_r$:
\begin{equation}\label{bounds}
\min_{B_r^d} H \leq \phi_d(x)\leq \max_{B_r^d} H, \; 0\leq \phi_i(x) \leq \max_{B_r^d} H-\min_{B_r^d}H,\; i=1,...,d-1.
\end{equation}
\end{prop}

\begin{proof}

\textbf{Step 1 : convex duality}

\smallskip

Equip $E:=C^{0}(B_r,\R)^d$ with the sup norm, and define (the linear continuous operator) $\Lambda$ by $\Lambda((\phi_1,...,\phi_d))(x_1,...,x_d):=\sum_{i=1}^d \phi_i(x_i)$ for all $(\phi_1,...,\phi_d)\in C^{0}(B_r,\R)^d$ and all $(x_1,...,x_d)\in B_r^d$. Remark now that \pref{dualhr} can be rewritten as:
\begin{equation}\label{duals}
\inf_{\phi\in E} F(\phi)+G(\Lambda(\phi))\end{equation}
with $F(\phi)= \sum_{i=1}^d \int_{B_r} \phi_i d\mu_i$ and, for all $\psi\in C^{0}(B_r^d,\R)$:
\begin{equation*}
G(\psi)=\left\{\begin{array}{lll}
0 &\mbox{ if }  \psi \geq H\\
+\infty  & \mbox{ otherwise }.
\end{array}
\right. \end{equation*}
The dual problem in the usual sense of convex analysis of \pref{duals} (see \cite{ektem}) is then
\begin{equation}\label{dduals}
\sup_{\gamma \in \MM(B_r^d)} -F^*(\Lambda^*\gamma)-G^*(-\gamma).\end{equation}
First, we remark that $\Lambda^*(\gamma)=(\pi_1\sharp \gamma,...,\pi_d\sharp \gamma)$. Elementary computations  then yield:
\begin{equation*}
F^*(\Lambda^*(\gamma))=\left\{\begin{array}{lll}
0 &\mbox{ if }  \pi_i\sharp \gamma=\mu_i \; \mbox{ for } i=1,...,d,\\
+\infty  & \mbox{ otherwise }.
\end{array}
\right. \end{equation*}
and 
\begin{equation*}
G^*(-\gamma)=\left\{\begin{array}{lll}
-\int_{\R^{d\times d}} H d\gamma &\mbox{ if }  \gamma \geq 0\\
+\infty  & \mbox{ otherwise }.
\end{array}
\right. \end{equation*}
Hence problem \pref{dduals} is  exactly \pref{mkh}. The Fenchel-Rockafellar duality theorem (see \cite{ektem})  implies then that  \pref{mkh}  admits solutions and that
\[\max \; \pref{mkh} = \inf \;  \pref{dualhr}.\]

\smallskip

\textbf{Step 2: convexification trick}

\smallskip

It remains to prove that the infimum is attained in \pref{dualhr}. Let $(\phi_1,...,\phi_d)\in \EE_r$ and define for all $x_1\in B_r$:
\begin{equation}\label{conj1}
\psi_1(x_1)=\sup_{(x_2,....,x_d)\in B_r^{d-1}} \left\{H(x_1,...,x_d)-\sum_{j=2}^d \phi_j(x_j) \right\} 
\end{equation}
by construction $\psi_1\leq \phi_1$ and $(\psi_1,\phi_2,...,\phi_d)\in \EE_r$. Construct then inductively $\psi_2,...,\psi_{d-1}$ by setting for $i=2,...,d-1$ and $x_i\in B_r$
\[\psi_i(x_i)= \sup_{(x_j)_{j\neq i} \in B_r^{d-1}} \left\{H(x_1,...,x_d)-\sum_{j<i} \psi_j(x_j)- \sum_{j>i} \phi_j(x_j)\right\}\]
and finally
\[\psi_d(x_d)= \sup_{(x_j)_{j\neq d} \in B_r^{d-1}} \left\{H(x_1,...,x_d)-\sum_{j=1}^{d-1} \psi_j(x_j) \right\}.\]
By construction, $\psi_i\leq \phi_i$ and $(\psi_1,...,\psi_d)\in\EE_r$ : $(\psi_1,...,\psi_d)$ is therefore an improvement of $(\phi_1,...,\phi_d)$ in problem \pref{dualhr}. On the one hand, since  $\psi_j\leq \phi_j$, one has for all $i=1,...,d$ and $x_i\in B_r$
\[\psi_i(x_i) \leq  \sup_{(x_j)_{j\neq i} \in B_r^{d-1}} \left\{H(x_1,...,x_d)-\sum_{j\neq i} \psi_j(x_j)\right\}.\]
On the other hand, the converse inequality holds because $(\psi_1,....,\psi_d)\in \EE_r$.

\smallskip

\textbf{Step 3: existence of minimizers for \pref{dualhr}}

\smallskip

Using the convexification trick of the previous step, we can find a minimizing sequence of \pref{dualhr}, $\phi^n:=(\phi_1^n,...,\phi_d^n)$ such that for all $n$, all $i$ and all $x_i\in B_r$, one has:
 \begin{equation}\label{conjn}
\phi_i^n(x_i)=\sup_{(x_j)_{j\neq i}\in B_r^{d-1}} \left\{H(x_1,...,x_d)-\sum_{j\neq i} \phi_j^n (x_j) \right\} 
\end{equation}
Noting that the objective in \pref{dualhr} is unchanged when changing $(\phi_1,...,\phi_d)$ into $(\phi_1+\alpha_1,...,\phi_d+\alpha_d)$ for constants $\alpha_i$ that sum to $0$, we may also assume that
\[\min_{B_r} \phi_i^n=0, \; \forall n\in \N,\; \forall i=1,...,d-1.\]
Together with \pref{conjn} we deduce that 
\[ \min_{B_r^d} H\leq  \phi_d^n \leq \max_{B_r^d} H \mbox{ on } B_r,  \; \forall n\in \N \]
and
\[ \phi_i^n \leq \max_{B_r^d} H- \min_{B_r^d} H.\]
Denoting by $\omega$ the modulus of continuity of $H$ on $B_r^d$, we also deduce  from \pref{conjn}, that for every $(x,y)\in B_r^2$, for every $n$ and $i$ one has:
\[\vert \phi_i^n(x)-\phi_i^n(y)\vert \leq \omega(\vert x-y\vert).\]
Thus, the sequence $\phi^n$ is bounded and uniformly equicontinuous hence by Ascoli's theorem admits some convergent subsequence. Denoting by $\phi=(\phi_1,...,\phi_d)$ the limit of this subsequence, it is easy to check that $\phi$ belongs to  $\EE_r$, solves \pref{dualhr} and satisfies \pref{conj} and \pref{bounds}.
\end{proof}


\subsection{The general case}

We now go back to $\MK$ (i.e. $H=\det$) for general marginals that only satisfy \pref{hypintegr}. By suitable truncation arguments, we have the following result, which is proved in Appendix B:

\begin{thm}\label{existd}
Assume that \pref{hypintegr} is satisfied, then both $\MK$ and $\DD$ admit solutions and
\[\max \; \MK = \min \; \DD.\]
Moreover, $\DD$ admits a solution $(\phi_1,....,\phi_d)$ such that  for all $i=1,...,d$ and all $x_i\in \R^d$, one has:
\begin{equation}\label{conjnc}
\phi_i(x_i)=\sup_{(x_j)_{j\neq i}\in \R^{d\times (d-1)}} \left\{\det(x_1,...,x_d)-\sum_{j\neq i} \phi_j(x_j) \right\}. 
\end{equation}

\end{thm}


\subsection{Extremality conditions}\label{cara}

This paragraph is devoted to optimality conditions for $\MK$ that can be derived from Theorem \ref{existd} (again we are in the case $H=\det$ here). For $(y_1,....,y_{d-1})\in (\R^d)^{d-1}$ we denote by $\bigwedge_{i=1}^{d-1} y_i$ the vector $y_1\wedge...\wedge y_{d-1}$ and recall that is characterized by the identity:
\[\det(y_1,....,y_d,x)= \<x,\bigwedge_{i=1}^{d-1} y_i>, \; \forall x\in\R^d.\]
For $(x_1,....,x_d)\in (\R^d)^{d}$ and $i\in\{1,...,d\}$, $\bigwedge_{j\neq i} x_j$ is defined in a similar way.

\smallskip

At this point, it is useful to remark that the family of l.s.c. functions $(\varphi_1,....,\varphi_d)$ belongs to $\EE$ if and only if for every $i\in\{1,...,d\}$ one has
\begin{equation}\label{caree}
\sum_{j\neq i} \varphi_j(x_j)\geq \varphi_i^*((-1)^{i+1} \bigwedge_{j\neq i} x_j),\; \forall (x_1,...,x_d)\in (\R^d)^d.
\end{equation}
It is also obvious that  $(\varphi_1,....,\varphi_d)$ belongs to $\EE$ if and only, it satisfies \pref{caree} for \emph{some} $i\in\{1,...,d\}$. Let us finally recall that the convexification trick ensures that one can always improve an element of $\EE$ in the dual problem by replacing it by another element of $\EE$ that satisfies \pref{conjnc}. Hence solutions of $\DD$ have to agree $\otimes_{i=1}^d \mu_i$ almost everywhere with potentials that satisfy \pref{conjnc}. Obviously if $(\varphi_1,....,\varphi_d)$  satisfies \pref{conjnc}, each $\varphi_i$ is a convex potential (as a supremum of a family of affine functions). With no loss of generality, we may therefore restrict ourselves to the subset $\EE_c$ consisting of the elements $(\varphi_1,....,\varphi_d)\in\EE$ such  that each potential $\varphi_i$ is convex on $\R^d$.

\smallskip

By definition of $\EE$ and the duality result of theorem \ref{existd}, we deduce that $\gamma\in \Pi(\mu_1,...,\mu_d)$ (respectively $(\varphi_1,...\varphi_d)\in \EE$) solves $\MK$ (respectively solves $\DD$) if and only if there exists $(\varphi_1,...\varphi_d)\in \EE$ (respectively $\gamma\in \Pi(\mu_1,...,\mu_d)$) such that
\begin{equation}\label{sellep}
\sum_{j=1}^d \varphi_j(x_j)=\det(x_1,...,x_d) \mbox{ $\gamma$-a.e.}
\end{equation} 
For all $\Phi:=(\varphi_1,...,\varphi_d)\in \EE$, let us define
\[K_{\Phi}:=\{(x_1,...,x_d)\in(\R^d)^d \mbox{ : }  
\sum_{j=1}^d \varphi_j(x_j)=\det(x_1,...,x_d)\}\]
and remark that  the (possibly empty) set $K_{\Phi}$ is closed since it is the minimal set of some l.s.c. function. Hence we deduce that \pref{sellep} is equivalent to $\gamma$ having its support included in $K_{\Phi}$. If $\Phi\in \EE_c$, we have the following characterization of $K_{\Phi}$:

\begin{lemma}\label{carkfi}
Let $\Phi:=(\varphi_1,...,\varphi_d)\in \EE_c$, and $x:=(x_1,...,x_d)\in (\R^d)^d$, then $x\in K_{\Phi}$ if and only if for all $i\in\{1,...,d\}$:
\begin{equation}\label{extr1}
\sum_{j\neq i} \varphi_j(x_j)=\varphi_i^*((-1)^{i+1} \bigwedge_{j\neq i} x_j)
\end{equation}
and 
\begin{equation}\label{extr2}
(-1)^{i+1} \bigwedge_{j\neq i} x_j \in \partial \varphi_i(x_i).
\end{equation}
\end{lemma}

\begin{proof}
If $x\in K_{\Phi}$ then
\[\sum_{j\neq i} \varphi_j(x_j)=\<x_i, (-1)^{i+1} \bigwedge_{j\neq i} x_j> -\varphi_i(x_i) \leq  \varphi_i^*((-1)^{i+1} \bigwedge_{j\neq i} x_j)\]
which proves \pref{extr1}. Now let $y_i\in\R^d$, since $x\in K_{\Phi}$ and $\Phi\in \EE_c$, we have:
\[\phi_i(y_i)-\phi_i(x_i)\geq \<y_i-x_i, (-1)^{i+1} \bigwedge_{j\neq i} x_j>,\]
which proves \pref{extr2}.

\smallskip

Conversely assume that $x$ satisfies \pref{extr1}-\pref{extr2} for some $i$. From \pref{extr1}, we get 
\[\sum_{k=1}^d \varphi_k(x_k)=\varphi_i^*((-1)^{i+1} \bigwedge_{j\neq i} x_j)+\varphi_i(x_i)\]
with \pref{extr2}, this yields
\[\sum_{k=1}^d \varphi_k(x_k)=\<x_i, (-1)^{i+1} \bigwedge_{j\neq i} x_j)>=\det(x),\]
so that $x\in K_{\Phi}$. 
\end{proof}
Remark that in fact,  for each fixed $i$, $K_{\Phi}$ consists of those $(x\in \R^{d})^{d}$ that satisfy \pref{extr1}-\pref{extr2} for that particular index $i$. We then immediately deduce the following:

\begin{prop}\label{fromprimaltodual}
Let $(\varphi_1,...,\varphi_d)\in \EE_c$, the following assertions are equivalent:

\begin{enumerate}

\item $(\varphi_1,...,\varphi_d)$ solves $\DD$,

\item there exists $\gamma\in \Pi(\mu_1,...,\mu_d)$ such that for every $i\in\{1,...,d\}$ and for $\gamma$-a.e. $(x_1,....,x_d)\in (\R^d)^d$, \pref{extr1} and \pref{extr2} hold,

\item  there exists $\gamma\in \Pi(\mu_1,...,\mu_d)$ and $i\in\{1,...,d\}$ such that for $\gamma$-a.e. $(x_1,....,x_d)$ \pref{extr1} and \pref{extr2} hold,

\end{enumerate}

\end{prop}

If the marginals have additional regularity, we immediately deduce the following uniqueness result for the dual problem:

\begin{prop}\label{uniq}
If for every $i\in\{1,...,d\}$, $\mu_i$ is absolutely continuous with respect to $\LL^d$  and has a positive density with respect to $\LL^{d}$ then $\DD$ admits a unique solution (up to the addition of constants summing to $0$ to each potential).
\end{prop}

\begin{proof}
Let ${\gamma}$ be a solution of $\MK$ and $\Phi:=(\phi_1,....,\phi_d)$ solve $\DD$, it follows from the duality relation that the support of $\gamma$ is included in $K_{\Phi}$. Moreover, it can be assumed that each $\varphi_i$ satisfies \pref{conjnc} hence $\Phi\in \EE_{c}$. Let $i\in\{1,...,d\}$, and let us desintegrate $\gamma$ with respect to its $i$-th marginal : $\gamma=\mu_i\otimes \gamma^{x_i}$. We deduce from lemma \ref{carkfi} that for almost every $x_i\in \R^d$ and $\gamma^{x_i}$ almost every $(x_j)_{j\neq i}\in (\R^d)^{d-1}$ one has:
\[ (-1)^{i+1} \bigwedge_{j\neq i} x_j \in \partial \varphi_i(x_i) \]
hence, by convexity
\[ (-1)^{i+1} \int_{(\R^d)^{d-1}} \bigwedge_{j\neq i} x_j  d\gamma^{x_i}((x_j)_{j\neq i}) \in \partial \varphi_i(x_i).\]
Since $\phi_i$ is differentiable $\mu_i$-a.e., we get that for $\mu_i$-a.e. $x_i$, one has 
\[\nabla \phi_i(x_i)=\int_{(\R^d)^{d-1}} \bigwedge_{j\neq i} x_j  d\gamma^{x_i}((x_j)_{j\neq i} )\]
since the rightmost member of this identity does not depend on $\Phi$, we are done.
\end{proof}

To sum up, getting back from $\DD$ to $\MK$, we obtain the following characterization of optimal transportation plans:

\begin{thm}\label{fromdualtoprimal}
Let $\gamma\in \Pi(\mu_1,...,\mu_d)$, $\gamma$ solves $\MK$ if and only if there exists l.s.c. convex functions $\varphi_i$ : $\R^d\to \R\cup \{+\infty\}$ such that for all $i\in\{1,...,d\}$:
\begin{eqnarray}
\sum_{j\neq i} \varphi_j(x_j)\geq \varphi_i^*((-1)^{i+1} \bigwedge_{j\neq i} x_j) \mbox{ on } (\R^d)^d,\\
\sum_{j\neq i} \varphi_j(x_j)= \varphi_i^*((-1)^{i+1} \bigwedge_{j\neq i} x_j) \mbox{ $\gamma$-a.e.},\\
(-1)^{i+1} \bigwedge_{j\neq i} x_j \in \partial \varphi_i(x_i) \mbox{ $\gamma$-a.e.}.
\end{eqnarray}
 
\end{thm}

Once again, one can replace ''for all $i\in\{1,...,d\}$'' by ''for some $i\in\{1,...,d\}$'' and "$\gamma$-a.e." by "on the support of $\gamma$'' in the previous result.  Of course, any $(\phi_1,...,\phi_d)$ satisfying the previous statements is a solution of $\DD$.

\smallskip

To illustrate the previous considerations, let us consider the case $d=3$ and assume that $(\varphi,\psi,\chi)$ is a (known)  triple of convex potentials that solve $\DD$. For the sake of simplicity, also assume that the marginals $(\mu_1,\mu_2,\mu_3)$ are absolutely continuous with respect to $\LL^{3}$. Then any optimal transport plan $\gamma$ is characterized by the extremality condition:
\[\varphi(x)+\psi(y)+\chi(z)=\det(x,y,z) \mbox{ on the support of $\gamma$}.\]
This implies  that, $\gamma$-a.e., one has
\[\left \{\begin{array}{lll}
\psi(y)+ \chi(z)=\varphi^*(y\wedge z)\\
\varphi(x)+\chi(z)=\psi^*(-x\wedge z)\\
\varphi(x)+\psi(y) =\chi^*(x\wedge y)
\end{array}\right.\]
and
\begin{equation}\label{supportopt}
\left \{\begin{array}{lll}
\nabla \varphi(x)=y\wedge z\\
\nabla \psi(y)=-x\wedge z\\
\nabla \chi(z)=x\wedge y
\end{array}\right.
\end{equation}
Note that in particular, one has:
\[\begin{split}
&\<x,\nabla \varphi(x)>=\<y,\nabla \psi(y)>=\<z,\nabla \chi (z) >\\
&=\det(x,y,z)=\det(\nabla \varphi(x),\nabla \psi(y),\nabla \chi(z)) \mbox{ $\gamma$-a.e. }.
\end{split}\] 
The  previous conditions clearly impose important geometric restrictions on $\gamma$.  It implies in particular that for $\mu_1$ almost every $x$, the conditional probability of $y$ or $z$ given $x$ is supported by $\nabla \phi(x)^{\perp}$. The conditional probability of $z$ given $(x,y)$ is even more constrained, indeed if $\<x,\nabla \varphi(x)>\neq 0$ then the previous conditions impose:
\[z=\frac{\nabla \varphi(x)\wedge \nabla \psi(y)}{\<x,\nabla \varphi(x)>}.\]

\section{The radial case}\label{radial}

In this section, we focus on $\MK$ in the case where the measures $(\mu_i)$ are radially symmetric. This means that for all $1\leq i\leq d$, and for all $R$ in the orthogonal group of $\R^d$,
$$
R\sharp\mu_i=\mu_i.
$$
Then, introducing the measures $\mu_i^r=|\cdot|\sharp\mu_i$, we have, for all $f\in C_c(\R^d)$,
\begin{equation}
\int_{\R^d}f(x)d\mu_i(x)=\frac{1}{|\S^{d-1}|}\int_{\S^{d-1}}\left(\int_0^{+\infty}f(re)d\mu_i^r(r)\right)d{\mathcal H}^{d-1}(e).
\end{equation}
In order to solve $\MK$, let us first remark that it is intuitive that the potentials which solve the dual problem $({\mathcal D})$ are radially symmetric. Then, noticing that \eqref{supportopt} generalizes to every dimension $d$, it follows that the support of an extremal measure will be included in the set of orthogonal systems. The only unknown here will be the relations between the norm of each vector. These relations will be obtained solving the following problem.
\begin{equation}
\tag{${\mathcal{MK}}^r$}
\label{symMK}
\sup_{\gamma^r\in\Pi(\mu_1^r,\ldots,\mu_d^r)}\left(\int_{(\R^d)^d}\left(\prod_{i=1}^dr_i\right)d\gamma^r(r_1,\ldots,r_d)\right).
\end{equation}
\begin{prop}
\label{symprob}
Assume that, for all $1\leq i\leq d$, $\mu_i^r$ has no atom. Then, the problem $({\mathcal{MK}}^r)$ admits a unique solution $\bar\gamma^r$ given by
$$
d\bar\gamma^r(r_1,\ldots,r_d)=d\mu_1^r(r_1)\otimes\delta_{\{r_2=H_2(r_1)\}}\otimes\ldots\otimes\delta_{\{r_d=H_d(r_1)\}},
$$
where, for each $2\leq i\leq d$, $H_i$ is the monotone rearrangement map of $\mu_1^r$ into $\mu_i^r$.
\end{prop}
\begin{proof}
Proceeding exactly as for $({\mathcal{MK}})$, we get the existence of a minimizer $\bar\gamma^r$ for $({\mathcal{MK}}^r)$, together with the existence of a solution for its dual problem
\begin{equation}
\label{symdual}
\tag{${\mathcal D}^r$}
\inf_{(\psi_1,\ldots,\psi_d)\in{\mathcal E}^r}\left(\sum_{i=1}^d\int_0^{+\infty}\psi_id\mu_i^r\right),
\end{equation}
where ${\mathcal E}^r$ is the set of $(\psi_1,\ldots,\psi_d)\in (C_0(\R_+))^d$, sucht that, for all $(r_1,\ldots,r_d)$ in $(\R_+)^d$,
$$
\sum_{i=1}^d\psi_i(r_i)\geq\prod_{i=1}^dr_i.
$$
As for $\MK$, if $(\psi_1,\ldots,\psi_d)$ is such a solution, we can assume that the $\psi_i$'s are l.s.c. convex functions on $\R^+$, and optimality conditions read as
\begin{itemize}
\item $\bar\gamma^r\in\Pi(\mu_1^r,\ldots,\mu_d^r)$.
\item for $\bar\gamma^r$-a.e. $(r_1,\ldots,r_d)\in\mbox{supp}(\bar\gamma^r)$,
\begin{equation}
\label{symopt}
\forall 1\leq i\leq d, \quad, \psi_i'(r_i)=\prod_{j\not=i}r_j.
\end{equation}
\end{itemize}
Let us introduce the functions $G_i(r)=r\psi_i'(r)$. Then, multiplying \eqref{symopt} by $r_i$, we get that, for $\bar\gamma^r$-a.e. $(r_1,\ldots,r_d)$ in $\mbox{supp}(\bar\gamma^r)$, and for all $1\leq i\leq d$,
$$
G_i(r_i)=G_1(r_1).
$$
Thanks to \eqref{symopt} again, it appears that each $\psi_i'$ is non negative and then, the function $G_i$ is non decreasing. It follows that, for $\mu_1^r$-a.e. $r_1$, and for all $(r_2,\ldots,r_d)$,
$$
\forall 1\leq i\leq d,\quad r_i=H_i(r_1):=(G_i^{-1}\circ G_1)(r_1),
$$
where $G_i^{-1}$ stands for the generalized inverse of $G_i$, and we get that the support of an optimal measure is necessarily in the closure of the graph of $(H_i)_{2\leq i\leq d}$, where all the $H_i$'s are non decreasing. But, the constraint on $\mbox{supp}(\bar\gamma^r)$ means exactly that, for all $2\leq i\leq d$, $H_i$ pushes $\mu_1^r$ forward to $\mu_i^r$. This ends the proof since, the measure $\mu_1^r$ being non atomic, such a map is unique.
\end{proof}
The main result of the section is the following.
\begin{thm}
\label{rad}
Let $(\mu_i)_{1\leq i\leq d}$ be radially symmetric measures on $\R^d$, and assume that the measures $\mu_i$ on $\R^+$ have no atom. Define the measure $\bar\gamma$ on $(\R^d)^d$ by
$$
\bar\gamma(x)=\mu_1(x_1)\otimes\bar\gamma^{x_1}(x_2)\otimes\bar\gamma^{x_1,x_2}(x_3)\otimes\ldots\otimes\bar\gamma^{x_1,\ldots,x_{d-1}}(x_d),
$$
with
\begin{equation}\label{coupling}
\left\{
\begin{array}{l}
\displaystyle\forall 1\leq i\leq d-2,\quad \bar\gamma^{x_1,\ldots,x_i}=\indi_{H_i(|x_1|)\S(x_1,\ldots,x_i)}\frac{d{\mathcal H}^{d-i-1}}{|\S^{d-i-1}|},\\
\displaystyle \S(x_1,\ldots,x_i)=\left\{x_1,\ldots,x_i\right\}^\perp\cap\S^{d-1},\\
\displaystyle\bar\gamma^{x_1,\ldots,x_{d-1}}=\delta_{\left\{H_d(|x_1|)\bigwedge_{i=1}^{d-1}\left(\frac{x_i}{|x_i|}\right)\right\}},
\end{array}
\right.
\end{equation}
with the maps $(H_i)_{2\leq i\leq d}$ given by Proposition \ref{symprob}. Then, $\bar\gamma$ is a solution to $({\mathcal{MK}})$.
\end{thm}
\begin{Remark} The previous solution is completely explicit since the $H_i$'s are.
\end{Remark}
\begin{Remark} As in section \ref{uniform}, the previous construction admits a very simple probabilistic interpretation. Indeed, the measure $\bar\gamma$ is the law of a vector $(X_1,\ldots,X_d)$ of random vectors, such that, for all $1\leq i\leq d$, $\mu_i$ is the law of $X_i$. Then, \eqref{coupling}  means that the conditional law of $X_i$ given $X_1=x_1,\ldots,X_{i-1}=x_{i-1}$ is uniform on $H(|x_1|)\S(x_1,\ldots,x_{i-1})$ for $1\leq i\leq d-1$, and $X_d$ is given by
$$
X_d=H_d(|X_1|)\bigwedge_{i=1}^{d-1}\left(\frac{X_i}{|X_i|}\right).
$$
\end{Remark}
\begin{proof}
It is sufficient to prove that $\bar\gamma\in\Pi(\mu_1,\ldots,\mu_d)$ and that there exists $(\varphi_1,\ldots,\varphi_d)\in{\mathcal E}$ such that \eqref{sellep} holds for all $(x_1,\ldots,x_d)$ in $\mbox{supp}(\bar\gamma)$. Notice that $\mbox{supp}(\bar\gamma)$ consists of the $(x_1,\ldots,x_d)$ such that
\begin{equation}
\label{radopt}
\left\{
\begin{array}{l}
\displaystyle (x_1,\ldots,x_d) \ \mbox{is an orthogonal basis of} \ \R^d,\\
\displaystyle \mbox{for} \ \mu_1 \ \mbox{a.e.} \ x_1, \forall 2\leq i\leq d,  |x_i|=H_i(|x_1|).
\end{array}
\right.
\end{equation}
Let us set
$$
\forall 1\leq i\leq d,\quad \varphi_i(x)=\psi_i(|x|),
$$
where $(\psi_1,\ldots,\psi_d)$ is a solution to the radial dual problem $({\mathcal D}^r)$. Then, for all $(x_1,\ldots,x_d)$ in $(\R^d)^d$,
$$
\sum_{i=1}^d\varphi(x_i)=\sum_{i=1}^d\psi(|x_i|)\geq\prod_{i=1}^d|x_i|\geq \mbox{det}(x_1,\ldots,x_d).
$$
In addition, If $(x_1,\ldots,x_d)$ belongs to the support of $\mbox{supp}(\bar\gamma)$, then by \eqref{radopt}, we both have
$$
\mbox{det}(x_1,\ldots,x_d)=\prod_{i=1}^d|x_i|
$$
and $(|x_i|)_{1\leq i\leq d}$ belongs to the support of the solution $\bar\gamma^r$ of $({\mathcal{MK}}^r)$ by proposition \ref{symprob}. Since $(\psi_i)_{1\leq i\leq d}$ is optimal in $({\mathcal D}^r)$, it follows that
$$
\mbox{det}(x_1,\ldots,x_d)=\sum_{i=1}^d\psi_i(|x_i|)=\sum_{i=1}^d\varphi_i(x_i).
$$
To end the proof, we just need to prove that $\bar\gamma$ has its marginals the $\mu_i$'s. The first marginal of $\bar\gamma$ is obviously $\mu_1$. Let $f\in C_0(\R^d)$, then by definition of $\bar\gamma$,
\begin{equation*}
\int_{\R^d}fd(\pi_2\sharp\bar\gamma)=\int_{\R^d}\left(\int_{H_2(|x_1|)\S(x_1)}f(x_2)\frac{d\HH^{d-2}(x_2)}{H_2(|x_1|)^{d-2}|\S^{d-2}|}\right)d\mu_1(x_1).
\end{equation*}
Using the radially symmetry of $\mu_1$, performing the change of variables $x_2=H_2(|x_1|)y$ in the inner integral, and then applying Lemma \ref{fubsphere}, we get
\begin{eqnarray*}
\int_{\R^d}fd(\pi_2\sharp\bar\gamma) & = & \int_0^\infty\left(\int_{\S^{d-1}}\left(\int_{\S(\sigma)}f(H_2(r)y)\frac{d\HH^{d-2}(y)}{|\S^{d-2}|}\right)\frac{d\HH^{d-1}(\sigma)}{|\S^{d-1}|}\right)d\mu_1^r(r)\\
& = & \int_0^\infty\left(\int_{\S^{d-1}}f(H_2(r)y)\left(\int_{\S(y)}\frac{d\HH^{d-2}(\sigma)}{|\S^{d-2}|}\right)\frac{d\HH^{d-1}(y)}{|\S^{d-1}|}\right)d\mu_1^r(r)\\
& = & \int_0^\infty\left(\int_{\S^{d-1}}f(H_2(r)y)\frac{d\HH^{d-1}(y)}{|\S^{d-1}|}\right)d\mu_1^r(r)\\
& = & \int_0^\infty\left(\int_{\S^{d-1}}f(\sigma)\frac{d\HH^{d-1}(\sigma)}{|\S^{d-1}|}\right)d(H_2\sharp\mu_1^r)(r)\\
& = & \int_{\R^d}fd\mu_2,
\end{eqnarray*}
by the definition of $H_2$. Repeated applications of this argument lead to the same result for the other marginals of $\bar\gamma$ up to the $(d-1)$-th. For the last one, we proceed as in the proof of Proposition \ref{3d-adm}. For a function $f\in C_0(\R^d)$, by change of variables,
\begin{eqnarray*}
\lefteqn{\int_{H_{d-1}(r)\S(x_1,\ldots,x_{d-2})}f\left(H_d(r)\bigwedge_{i=1}^{d-1}\left(\frac{x_i}{|x_i|}\right)\right)\frac{d\HH^1(x_{d-1})}{H_d(r)|\S^1|}} & & \\
& & =\int_{\S(x_1,\ldots,x_{d-2})}f\left(H_d(r)\bigwedge_{i=1}^{d-1}\left(\frac{x_i}{|x_i|}\right)\right)\frac{d\HH^1(x_{d-1})}{|\S^1|},\\
& & =\int_{\S(x_1,\ldots,x_{d-2})}f\left(H_d(r)x_d\right)\frac{d\HH^1(x_d)}{|\S^1|},
\end{eqnarray*}
noticing that, if $x_1,\ldots,x_{d-2}$ are fixed, the map
$$
x_{d-1}\mapsto\bigwedge_{i=1}^{d-1}\left(\frac{x_i}{|x_i|}\right)
$$
is a rotation with angle $\pi/2$ on $\S(x_1,\ldots,x_{d-2})$. We can then repeat the argument used for the other marginals and then finish the proof.
\end{proof}

\section{Volume maximization}

So far, we have restricted our attention to the case where the objective function is the determinant although we were initially motivated with a volume maximization problem which corresponds to:
\[\MK_a \;  \sup_{\gamma \in \Pi(\mu_1,...,\mu_d)} \int_{\R^{d\times d}} \vert \det (x_1,...,x_d)\vert d\gamma(x_1,...,x_d).\]
instead of $\MK$. In the radial case, treated in the previous section, we have seen that optimal measures for $\MK$ give full mass to the set of matrices with nonnegative determinant. In this case, there is no loss of generality in replacing $\MK_a$ by $\MK$ in the sense that solutions of $\MK$ also solve $\MK_a$. This holds true under less restrictive symmetry assumptions on the marginals:

\begin{prop}
Assume that \pref{hypintegr} is satisfied and that  among the marginals $(\mu_1,..., \mu_d)$, two are symmetric (i.e. $(-{\rm{id}})\sharp \mu_i=\mu_i$, $(-{\rm{id}})\sharp \mu_j=\mu_j$ for two different indices $i$ and $j$) then any solution of $\ogm$ of $\MK$ satisfies $\det(x_1,...,x_d)\geq 0-\ogm$-a.e. and actually also solves $\MK_a$.
\end{prop}

\begin{proof}
Without loss of generality, assume that $\mu_1$ and $\mu_2$ are symmetric and that $d\geq 3$. Assume that $(\varphi_1,....,\varphi_d)\in\EE_c$ solves the dual problem $\DD$, then so does $(\varphi_1\circ(-\rm{id}), \varphi_2\circ(-\rm{id}),\varphi_3,...,\varphi_d)$. Hence $(\psi_1,...,\psi_d)$ also solves $\DD$ where $\psi_1$ and $\psi_2$ are the even part of $\varphi_1$ and $\varphi_2$ respectively and $\psi_i=\varphi_i$ for $i\geq 3$. Since $\psi_1$ is even we actually have for all $(x_1,...,x_d)$:
\begin{equation}\label{valabsd}
\psi_1(-x_1)+\sum_{j=2}^d \psi_j(x_j)=\sum_{i=1}^d \psi_i(x_i)\geq \vert \det(x_1,...,x_d)\vert \end{equation}
but, from the duality relations, one deduces that $\ogm$-a.e.
\[\sum_{i=1}^d \psi_i(x_i)= \det(x_1,...,x_d)\geq 0.\]
To show that $\ogm$ solves $\MK_a$, we use a classical duality argument. Indeed, let $\gamma\in\Pi(\mu_1,...,\mu_d)$, from \pref{valabsd}, one deduces:
\[\begin{split}
\int_{\R^{d\times d}} \vert \det (x_1,...,x_d)\vert d\gamma(x_1,...,x_d)\leq \sum_{i=1}^d \psi_i(x_i) d\mu_i(x_i)\\
=\int_{\R^{d\times d}} \vert \det (x_1,...,x_d)\vert d\ogm(x_1,...,x_d).
\end{split}\]

\end{proof}

\section{Uniqueness issues}\label{comments}

We end the paper by some remarks on uniqueness. For this, we will mainly use the very simple example of section \ref{uniform}. First, as noticed at the end of section \ref{uniform}, if the objective function is
$$
x\mapsto \vert\det(x)\vert,
$$
then there are  infinitely many possible choices for the third random vector, given the  first $2$ ones. Indeed, for all $x$ and $y$ in the ball $B$ with $\vert x\vert =\vert y \vert$, and for all $p\in[0,1]$, then we can choose for the conditional law of $Z$, given $X=x$ and $Y=y$:
$$
p\delta_{\frac{x\wedge y}{\vert x\vert}}+ (1-p)\delta_{\frac{y\wedge x}{\vert x\vert}}.
$$
Note that for the determinant, only the value $p=1$ is optimal and it corresponds to the solution given by Theorem \ref{thm-3d}. In the sequel, we shall work with the determinant as objective function.

\medskip

As seen in the introduction, in dimension $2$, since the problem $\MK$ can be reduced to the classical Monge-Kantorovich problem with quadratic cost, the solution is unique, at least for sufficiently regular measures. Actually, we have more, that is, up to a rotation, the optimal measure is supported by the graph of the gradient of a convex function $\varphi$. In terms of random pairs $(X,Y)$, the optimal coupling is characterized by the relation $RY=\nabla\varphi(X)$, with $\varphi$ convex and $R$ the rotation with angle $-\pi/2$. In higher dimensions,  such uniqueness is lost.

\medskip

First, let us consider the solution given by theorem \ref{thm-3d}. For $x\in B$, and given $X=x$, we define the conditional law of $Y$ by the measure on $\vert x\vert\S(x)$, absolutey continuous with respect to the $1$-dimensional Hausdorf measure, with density
$$
\rho(y)=1+\frac{\< y,e>}{|y|},
$$
where $e$ is an arbitrary fixed vector in $\S^2$, which is a probability measure since we add to $1$ a function which is greater than $-1$ and has zero mean on any circle $\S(x)$. The vector $Z$ is then chosen exactly as in Theorem \ref{thm-3d}.
The proof of the optimality of this measure is the same as for Theorem \ref{thm-3d}.
\medskip

Finally, we study the possibility for the extremal measure to have the same structure as in the 2-dimensional case, that is to be supported by the graph of a map $T:\R^d\to(\R^d)^{d-1}$. Such a $T$ would then be called a ''Monge'' solution. Notice that the solution given by Theorems \ref{thm-3d} and \ref{rad} does not have this structure but the question of the existence of such a solution is quite natural. Let us set the problem on $\R^d$ with all the measures uniform on the ball, which  is actually covered by Theorem \ref{rad}. Then, according to the extremality conditions, the existence of a Monge solution for $\MK$ is equivalent to the existence of $(d-1)$ maps $T_2,\ldots T_{d-1}$: $B \to B$ such that:
\begin{equation}\label{monge}
 \left\{
 \begin{array}{l}
 \mbox{for a.e. }  x\in B, \ \left\{x,T_2(x),\ldots,T_d(x)\right\} \ \mbox{is an orthogonal basis of} \ \R^d, \\
 \mbox{for a.e. }  x\in B, \forall 2\leq i\leq d, \ |T_i(x)|=|x|,\\
  \forall 2\leq i\leq d, T_i \ \mbox{preserves the Lebesgue measure}.
 \end{array}
\right.
\end{equation}
The construction of such maps seems to be very difficult  in general and we left it open. However, in dimension $4$, we have the following explicit construction:
\begin{equation*}
 T_2(x)=\left(
 \begin{array}{c} -x_2 \\ x_1 \\ -x_4 \\ x_3
 \end{array}
 \right), \
 T_3(x)=\left(
 \begin{array}{c} -x_3 \\ x_4 \\ x_1 \\ -x_2
 \end{array}
 \right), \
 T_4(x)=\left(
 \begin{array}{c} -x_4 \\ -x_3 \\ x_2 \\ x_1
 \end{array}
 \right).
 \end{equation*}
In contrast, in dimension $3$, it is not clear to us whether one can even find a single measure preserving map $T$ : $B\rightarrow B$ such $T$ is norm-preserving and a.e. orthogonal  to the identity map. If such a map exists, then by the hairy ball Theorem, it cannot be continuous.

\section*{Appendix A}

\begin{lemma}\label{fubsphere}
Let $k\geq 2$ be an integer and $f\in C_0(\S^k\times\S^k)$. Then,
\begin{equation}
\label{fubini}
\int_{\S^k}\left(\int_{\S(x)}f(x,y)\frac{d\HH^{k-1}(y)}{|\S^{k-1}|}\right)\frac{d\HH^k(x)}{|\S^k|} = \int_{\S^k}\left(\int_{\S(y)}f(x,y)\frac{d\HH^{k-1}(x)}{|\S^{k-1}|}\right)\frac{d\HH^k(y)}{|\S^k|},
\end{equation}
where
$$
\S(\xi)=\left\{ z\in\S^k ; \< z,\xi>=0\right\}.
$$
\end{lemma}

\begin{proof}
Let us introduce the following notations. Let us define the measure $\mu$ on $\S^k\times\S^k$ by
$$
d\mu(x,y)=d\nu(x)\otimes d\theta^x(y),
$$
where
$$
d\nu(x)=\frac{d\HH^k(x)}{|\S^k|} \ \mbox{and} \ d\theta^x(y)=\indi_{\S(x)}(y)\frac{d\HH^{k-1}(y)}{|\S^{k-1}|}.
$$
The identity \eqref{fubini} expresses the fact that $\mu$ is symmetric in $x$ and $y$. We set, for $x\in\S^k$ and $\varepsilon>0$,
$$
\S_\varepsilon(x)=\left\{ z\in\S^k; \ |\< z,x>|\leq\varepsilon \right\}.
$$
Let us notice first that $\HH^k(\S_\varepsilon(x))$ is independant of $x$ and denote this quantity by $\omega_\varepsilon$. Indeed, for all $x$ and $x'$ in $\S^k$, let $R$ be a rotation that maps $x$ to $x'$. Then, we have
$$
\S_\varepsilon(x')=\S_\varepsilon(Rx)=R^{-1}\left(\S_\varepsilon(x)\right),
$$
and then $\HH^k(\S_\varepsilon(x'))=\HH^k(\S_\varepsilon(x))$ since the Hausdorff measure is invariant under rotations. Next, we define the probability measure $\theta^x_\varepsilon$ on $\S^k$ as
$$
d\theta^x_\varepsilon(y)=\indi_{\S_\varepsilon(x)}(y)\frac{d\HH^k(y)}{\omega_\varepsilon}.
$$
For all $f\in C_0(\S^k,\S^k)$, and for all $x\in\S^k$,
$$
\left|\int_{\S^k}f(x,y)d\theta^x_\varepsilon(y)-\int_{S^k}f(x,y)d\theta^x(y)\right| \leq \left|\int_{\S^k}\left(f(x,y)-f(x,T^x(y))\right)d\theta^x_\varepsilon(y)\right|
$$
\begin{equation}
\label{decomp} + \left|\int_{\S^k}f(x,T^x(y))d\theta^x_\varepsilon(y)-\int_{S^k}f(x,y)d\theta^x(y)\right|,
\end{equation}
where
$$
T^x(y)=\frac{\pi_xy}{\|\pi_xy\|},
$$
and $\pi_x$ is the orthogonal projection onto $x^\perp$. First, we shall estimate the first term in the right-hand side of \eqref{decomp}. It is not difficult to see that if $\varepsilon\leq 1/2$, then for all $x\in\S^k$ and $y\in\S_\varepsilon(x)$,
$$
\|y-T^x(y)\|\leq 4\varepsilon,
$$
thus, denoting by $\omega_f$ the modulus of continuity of $f$,
$$
\left|\int_{\S^k}\left(f(x,y)-f(x,T^x(y))\right)d\theta^x_\varepsilon(y)\right|\leq \omega_f(4\varepsilon).
$$
For the second term in the right-hand side of \eqref{decomp}, we write
$$
\int_{\S^k}f(x,T^x(y))d\theta^x_\varepsilon(y) = \int_{\S^k}f(x,z)d(T^x\sharp\theta^x_\varepsilon)(z),
$$
where
$$
T^x(y)=\frac{\pi_xy}{\|\pi_xy\|}.
$$
Since $T^x$ commutes with any rotation that leaves $x$ invariant and $\theta^x_\varepsilon$ is invariant under such a rotation, the measure $T^x\sharp\theta^x_\varepsilon$ has its support in $\S(x)$ and is invariant under rotation on $\S(x)$. So necessarily, we have
$$
T^x\sharp\theta^x_\varepsilon=\theta^x
$$
and
$$
\left|\int_{\S^k}f(x,y)d\theta^x_\varepsilon(y)-\int_{S^k}f(x,y)d\theta^x(y)\right| \leq \omega_f(4\varepsilon).
$$
As a consequence,
\begin{equation}
\label{cv-unif}
\lim_{\varepsilon\to 0}\left(\sup_{x\in\S^k}\left|\int_{\S^k}f(x,y)d\theta^x_\varepsilon(y)-\int_{S^k}f(x,y)d\theta^x(y)\right|\right)=0.
\end{equation}
We set $d\mu_\varepsilon(x,y)=d\nu(x)\otimes d\theta^x_\varepsilon(y)$. Then, \eqref{cv-unif} implies that $(\mu_\varepsilon)_\varepsilon$ converges weakly to $\mu$. In addition, $\mu_\varepsilon$ is symmetric in $x$ and $y$, which ends the proof.
\end{proof}

\section*{Appendix B}

\begin{proof}[ of theorem \ref{existd}]

\textbf{Step1 : existence of maximizers for $\MK$}

\smallskip

Let $k\in\N$ and define:
\begin{equation}\label{defundHk}
\undH_k:=\max(\undH,-k)\end{equation}
By construction, $-k\leq \undH_k\leq 0$ and $\undH_k$ decreases to $\undH$.

Let $\gamma_n$ be a maximizing sequence of $\MK$. Since $\Pi(\mu_1,...,\mu_d)$ is tight, it follows from Prohorov's Theorem, that (taking a subsequence if necessary) there exists $\gamma\in \MM(\R^{d\times d})$ such that for every continuous bounded function $f$ on $\R^{d\times d}$,  $\idd f \gamma_n$ converges to $\idd f d\gamma$. Clearly, $\gamma\in \Pi(\mu_1,...,\mu_d)$ and for every $k$ 
\[\idd \undH_k d\gamma=\lim_n \idd \undH_k d\gamma_n \geq \lim_n \idd \undH d\gamma_n=\sup \MK.\]
It thus follows from  the monotone convergence theorem that $\gamma$ is a solution of  $\MK$.

\smallskip

\textbf{Step 2 : duality by approximation}

\smallskip

Consider the problem:
\begin{equation*}
\MKk \;  \sup_{\gamma \in \Pi(\mu_1,...,\mu_d)} \int_{\R^{d\times d}} \undH_k(x) d\gamma(x)
\end{equation*}
By the same considerations as in step 1, $\MKk$ admit solutions, let $\gamma_k$ be such a solution. Now for $\eps>$, let $r=r_\eps>0$ be such that 
\begin{equation}\label{tightre}
\sum_{i=1}^d \mu_i(\R^d \setminus B_r^d)+ \idpbd (1+\vert H_0(x)\vert)d\gamma(x)\leq \eps, \; \forall \gamma\in \Pi(\mu_1,...,\mu_d).
\end{equation}
Define then for $i=1,...,d$:
\[\mu_i^{r,k}:=\pi_i \sharp\left ( \frac { \gamma_k \indi_{B_r^d}  }{ \gamma_k(B_r^d)} \right) .\]  
Consider now:
\begin{equation*}
\MKrk \;  \sup_{\gamma \in \Pi(\mu_1^{r,k},...,\mu_d^{r,k})} \int_{\R^{d\times d}} \undH_k(x) d\gamma(x)
\end{equation*} 
and (defining $\EE_{r,k}$ as in \pref{contraintedudualr} with $\undH_k$ instead of $\det$) its dual:
\begin{equation*}
\DDrk \;  \inf_{(\phi_1,...,\phi_d)\in\EE_{r,k}}\sum_{i=1}^d \int_{\R^d} \phi_i d\mu_i^{r,k}. \end{equation*}
We know from proposition \ref{cascompact} that the values of $\MKrk$ and $\DDrk$ are equal and both attained. Moreover, using again proposition \ref{cascompact}, the infimum in $\DDrk$ is attained by some $\phi^{r,k}=(\phi_1^{r,k}, ..., \phi_d^{r,k})$ such that:
\begin{equation}\label{conjrk}
\phi_i^{r,k}(x_i)=\sup_{(x_j)_{j\neq i}\in B_r^{d-1}} \left\{\undH_k(x_1,...,x_d)-\sum_{j\neq i}\phi_j^{r,k}(x_j) \right\},\; \forall x_i\in B_r, 
\end{equation}
and, for all $x\in B_r$:
\begin{equation}\label{boundsrk}
-k \leq \phi_1^{r,k}(x)\leq 0, \; 0\leq \phi_i^{r,k}(x) \leq k ,\; i=2,...,d.
\end{equation}
Define then for all $x_1\in\R^d$:
\begin{equation}\label{conj1r}
\psi_1^{r,k}(x_1)=\sup_{(x_2,....,x_d)\in B_r^{d-1}} \left\{\undH_k(x_1,...,x_d)-\sum_{j=2}^d \phi^{r,k}_j(x_j) \right\} 
\end{equation}
Construct then inductively $\psi_2^{r,k},...,\psi_{d-1}^{r,k}$ by setting for $i=2,...,d-1$ and $x_i\in \R^d$,
$$
\psi_i^{r,k}(x_i)= \sup \{\undH_k(x_1,...,x_d)-\sum_{j<i} \psi_j^{r,k}(x_j)- \sum_{j>i} \phi_j^{r,k}(x_j),
$$
$$
(x_1,...,x_{i-1},x_{i+1},...,x_d)\in \R^{d\times (i-1)}\times B_r^{d-i} \}
$$
and finally
\[\psi_d^{r,k}(x_d)= \sup_{(x_j)_{j\neq d} \in \R^{d\times(d-1)}} \left\{\undH_k(x_1,...,x_d)-\sum_{j=1}^{d-1} \psi_j^{r,k}(x_j) \right\}.\]
It can be checked easily that each $\psi_i^{r,k}$ extends  $\phi_i^{r,k}$ to the whole of $\R^d$. By \pref{boundsrk} it can also be checked that $\Vert \psi^{r,k} \Vert_{\infty}=O(k)$. By construction $\psi^{r,k}=(\psi_1^{r,k}, ..., \phi_d^{r,k})\in \EE^k$, where $\EE_k$ is the set of $d$-uples of lower-semi continuous functions $(\phi_1,...,\phi_d)$ from $\R^d$ to $\R\cup \{+\infty\}$ such that
\begin{equation}\label{contraintedudualk}
 \sum_{i=1}^d \phi_i(x_i)\geq \undH_k(x_1,...,x_d),\; \forall \; (x_1,...,x_d)\in\R^{d\times d}. 
\end{equation}
In other words, $\EE^k$ is the admissible set of the dual problem of $\MKk$:
\begin{equation*}
\DDk \;  \inf_{(\phi_1,...,\phi_d)\in\EE^{k}}\sum_{i=1}^d \int_{\R^d} \phi_i d\mu_i. \end{equation*}
Let $\gamma_{r,k}$ be a solution of $\MKrk$ and define:
\[\ovga_{r,k}:=\gamma_k(B_r^d) \gamma_{r,k}+ \indi_{\R^{d\times d}\setminus B_r^d}.\]
It is easy to verify that $\ovga_{r,k} \in \Pi(\mu_1,...,\mu_d)$ hence we get:
\[\begin{split} 
\sup \MKk  =  &\idd \undH_k d\gamma_k \geq  \idd \undH_k d\ovga_{r,k} \\
=  & \gamma_k(B_r^d) \int_{B_r^d} \undH_k d\gamma_{r,k}+\idpbd \undH_k d\gamma_k\\
= & \gamma_k(B_r^d)\left( \sum_{i=1}^d \int_{B_r} \phi_i^{r,k} d\mu_i^{r,k}  \right) +\idpbd \undH_k d\gamma_k
\end{split}\]
 Using \pref{tightre}, we thus get:
 \[\begin{split}  \sup \MKk  \geq    \gamma_k(B_r^d)\left( \sum_{i=1}^d \int_{B_r} \phi_i^{r,k} d\mu_i^{r,k}  \right)  -\eps=   \sum_{i=1}^d \int_{B_r^d} \phi_i^{r,k} d\gamma_k  -\eps\\
 =  \sum_{i=1}^d \idd  \psi_i^{r,k} d\gamma_k  -  \sum_{i=1}^d \idpbd \psi_i^{r,k} d\gamma_k    -\eps \geq  \inf  \DDk +O(k) \eps. 
\end{split}\]
Letting $\eps$ go to $0$ we then get $\sup \MKk\geq \inf  \DDk$. Using Prohorov's Theorem, taking a subsequence if necessary, we may assume that  there exists $\gamma\in \Pi(\mu_1,..,\mu_d)$ such that for every continuous bounded function $f$ on $\R^{d\times d}$,  $\idd f \gamma_k$ converges to $\idd f d\gamma$. Since $\undH_k\geq \undH=\det-H_0$, we then have for all $k$:
\[\sup \MKk=\idd \undH_k d\gamma_k \geq \inf \DD-a\]
Since $\undH_k$ is nonincreasing in $k$, for all $k_0$, we then have:
\[\idd \undH_{k_0} d\gamma =\lim_k \idd  \undH_{k_0} d\gamma_k \geq  \idd \undH_k d\gamma_k \geq \inf \DD-a\]
taking the infimum in $k_0$ and using the monotone convergence theorem we then get:
\[ \sup\MK -a \geq \idd \undH d\gamma\geq \inf \DD-a.\]
Since we already know that $\sup\MK \leq \inf \DD$, we have $\sup\MK =\inf \DD$.

\smallskip

\textbf{Step  3: existence of minimizers for $\DD$}

\smallskip

It remains to prove that the infimum is attained in $\DD$. By the convexification trick,  we can find a minimizing sequence of $\DDr$, $\phi^n:=(\phi_1^n,...,\phi_d^n)$ such that for all $n$, all $i$ and all $x_i\in \R^d$, one has:
 \begin{equation}\label{conjnm}
\phi_i^n(x_i)=\sup_{(x_j)_{j\neq i}\in \R^{d\times{d-1}}} \left\{\det (x_1,...,x_d)-\sum_{j\neq i} \phi_j^n (x_j) \right\} 
\end{equation}
Note the functions $\phi_i^n$ are convex l.s.c. and not identically equal to $+\infty$. By \pref{conjnm}, $\phi_i^n$ admits therefore some affine minorant hence $\phi_i^n(x)+\vert x \vert^{p_i}/p_i$ achieves its minimum. Using again  that the objective in $\DD$ is unchanged when changing $(\phi_1,...,\phi_d)$ into $(\phi_1+\alpha_1,...,\phi_d+\alpha_d)$ for constants $\alpha_i$ that sum to $0$, we may then also assume that for all $i=1,...,d-1$:
\begin{equation}\label{normali}
\min_{x\in \R^d} \left\{ \phi_i^n(x)+\frac{ \vert x \vert^{p_i}}{p_i}\right\}=0.
\end{equation}
We easily deduce from the previous and \pref{conjnm} that 
\[\phi_d^n(x)+\frac{ \vert x \vert^{p_d}}{p_d}\geq 0, \; \forall x\in \R^d.\]
Let $m>0$ and define truncated potentials $\phi_i^{n,m}$ and $\psi_i^{n,m}$ by:
\begin{equation}\label{tronc}
\psi_i^{n,m}(x)=\phi_i^{n,m}(x)+\frac{ \vert x \vert^{p_i}}{p_i}=\min(\phi_i^{n}(x)+\frac{ \vert x \vert^{p_i}}{p_i},m)\; \forall x\in \R^d.
\end{equation}
By construction, $0\leq \psi_i^{n,m} \leq m$ and there exists $C$ such that for all $(n,m)$, one has:
\begin{equation}\label{borneint}
\sum_{i=1}^{d} \idd \psi_i^{n,m} d\mu_i \leq C.
\end{equation}
and since $\psi_i^{n,m}\geq 0$, for all $(x_1,...,x_d)\in \R^{d\times d}$, one has:
\begin{equation}\label{somsupnm}
\sum_{i=1}^{d}  \psi_i^{n,m} (x_i)\geq \min(\ovH(x_1,...,x_d), m).
\end{equation}
For fixed $m$, taking subsequence if necessary, we may assume that $\psi_i^{n,m}$ weakly converges in $L^{1}(\mu_i)$ to some limit $\psi_i^m$ as $n\rightarrow+\infty$. By Mazur's Lemma, there is some sequence of convex combinations of the  $(\psi_i^{l,m})_{l\geq n}$ that  converges in $L^{1}(\mu_i)$ and (possibly after an extraction) $\mu_i$-almost everywhere to $\psi_i^m$ as $n\rightarrow+\infty$, there is no loss of generality of setting $\psi_i^m=+\infty$ outside the set where there is convergence. Hence, passing to the limit in \pref{somsupnm}, we get, for all $(x_1,...,x_d)\in \R^{d\times d}$
\begin{equation}\label{somsupm}
\sum_{i=1}^{d}  \psi_i^{m} (x_i)\geq \min(\ovH(x_1,...,x_d), m).
\end{equation}
Now, $\psi_i^m$ is a nondecreasing sequence of nonnegative functions that is bounded in $L^{1}(\mu_i)$. By the monotone convergence theorem, $\psi_i^m$ converges pointwise and in $L^1(\mu_i)$ to some function $\psi_i$. Defining $\phi_i(x):=\psi_i(x)-\vert x \vert^{p_i}/p_i$, and passing to the limit in \pref{somsupm} first yields:
\begin{equation}\label{somsuplim}
\sum_{i=1}^{d}  \phi_i^{m} (x_i)\geq \det(x_1,...,x_d), \; \forall (x_1,...,x_d)\in \R^{d\times d}.
\end{equation}
Since $\psi_i^n\geq \psi_i^{n,m}$, we also get:
\[\begin{split}
\inf \DD = &\lim_n \sum_{i=1}^{d}\idd  \phi_i^{n}d\mu_i \geq \sup_m \left ( \lim_n   \sum_{i=1}^{d}\idd \psi_i^{n,m}d\mu_i \right) -a \\
= & \sum_{i=1}^{d}\idd  \psi_i d\mu_i-a =\sum_{i=1}^{d}\idd  \phi_id\mu_i 
\end{split}\]
Thanks to \pref{somsuplim}, we can apply the convexification trick to $(\phi_1,...,\phi_d)$. We then obtain an element of $\EE$ that satisfies \pref{conjnm} and solves $\DD$.

\end{proof}

\end{document}